\documentclass[a4paper,10pt]{article}

 \usepackage{a4wide}

 \usepackage{graphics,graphicx}
 \usepackage{amsmath,amsthm}
 \usepackage{amsfonts,mathrsfs,amssymb}
 \newtheorem{theorem}{Theorem}[section]
 \newtheorem{lemma}[theorem]{Lemma}
 \newtheorem{proposition}[theorem]{Proposition}
 \newtheorem{definition}[theorem]{Definition}

 \theoremstyle{definition}
 \newtheorem{remark}[theorem]{Remark}

 \newcommand{\pr}{\noindent{\bf Proof. }}
 \newcommand{\ep}{\nolinebreak{\hspace*{\fill}$\Box$ \vspace*{0.25cm}}}

 \newcommand{\mb}[1]{\ensuremath{\mathbb{#1}}}
 \newcommand{\N}{\mb{N}}
 
 \newcommand{\R}{\mb{R}}
 \newcommand{\C}{\mb{C}}
 
 \newcommand{\LT}[1]{\ensuremath{{\cal L}#1}}
 \newcommand{\ILT}[1]{\ensuremath{{\cal L}^{-1}#1}}

 \newcommand{\al}{\alpha}
 \newcommand{\be}{\beta}
 \newcommand{\ga}{\gamma}
 \newcommand{\del}{\delta}
 \newcommand{\Ga}{\Gamma}

 \newcommand{\wh}{\widehat}
 
 \newcommand{\vphi}{\varphi}

 \newcommand{\spp}{\mathscr{S}'_+}
 \newcommand{\dpp}{\mathscr{D}'_+}
 
 \newcommand{\E}{\ensuremath{{\cal E}}}
 
 \renewcommand{\S}{\mathscr{S}}
 \newcommand{\Lloc}{L_{\text{loc+}}^{1}}
 \newcommand{\supp}{\mathop{\mathrm{supp}}}
 \renewcommand{\Re}{{\rm {Re\,}}}

 \newcommand{\dis}[2]{\langle #1 , #2 \rangle}

 \title{Semilinear ordinary differential equation coupled
 with distributed order fractional differential equation}

 \author{Teodor M. Atanackovic\thanks{Faculty of Technical Sciences, University of Novi Sad,
  Trg D. Obradovi\'ca 5, 21000 Novi Sad, Serbia,
  (\tt{atanackovic@uns.ns.ac.yu})}\hspace{0,2cm}
  Ljubica Oparnica
  \thanks{Institute of Mathematics, Serbian Academy of Science,
  Kneza Mihaila 35, 11000 Belgrade, Serbia,
  (\tt{ljubicans@sbb.co.yu})}\hspace{0,2cm}
  Stevan Pilipovi\'{c}
  \thanks{Department of Mathematics and Informatics, University of Novi Sad,
  Trg Dositeja Obradovi\'ca 4, 21000 Novi Sad, Serbia,
  (\tt {pilipovic@im.ns.ac.yu})}
  }

\begin{document}
 \date{}
 \maketitle

 \begin{abstract}
 \noindent System
 $D^{2}y(t)+z(t)=f(t,y)$,
 $\int_{0}^{2}\phi _{1}(\gamma)D^{\gamma}y(t)d\gamma =
 \int_{0}^{2}\phi _{2}(\gamma )D^{\gamma}z(t)d\gamma$, $t>0, $
 where $D^{\gamma}, \gamma\in [0,2]$ are  fractional derivatives,
 is investigated and the existence of the solution
 in a mild and classical sense is proven.
 Such a system arises as
 a distributed derivative model of viscoelastic body and
 in the system identification theory.
 Also, the condition for existence and uniqueness
 of a solution to a general linear fractional differential equation
 $\sum_{i=1}^k a_iD^{\ga_i}z(t)=g(t)$,  $\ga_i\in\R,\,i=1,...,k$
 in $ \spp$  is given.\\

% \begin{keywords}
 %\end{keywords}

 %\begin{AMS}
 \noindent {\footnotesize AMS Subject
 Classification Subject classification: 26A33; 34G20; 47H10}
 %\end{AMS}
 \end{abstract}

 \pagestyle{myheadings} \thispagestyle{plain}
 \markboth{T.M.Atanackovi\'c, Lj. Oparnica and
 S.Pilipovi\'c}{System with distributed fractional derivatives}

 \section{Introduction}\label{intro}

  %\subsection{Notation}

 In this paper, we  solve and  analyze solutions to a nonlinear
 system motivated by a mechanical model
 \begin{equation*}
 (CP)\quad\quad
 \left\{ \begin{array}{c}
  D^{2}y(t)+z(t)=f(t,y(t)),\quad t > 0, \\
   \\
  \int_{0}^{2}\phi _{1}(\gamma )D^{\gamma }y(t)d\gamma =
   \int_{0}^{2}\phi_{2}(\gamma )D^{\gamma }z(t)d\gamma,\quad t > 0, \\
    \\
  y(0)=y_{0};\qquad y^{(1)}(0)=v_{0}.
  \end{array}
 \right.
 \end{equation*}
  Here, $y$ and $z$, which represent strain and stress respectively,
  are unknown functions defined for $t>0$. Given are locally integrable function
 $f$ representing  forcing term,  and $\phi _{1},\phi _{2}$
 %which are certain functions or distributions
 which characterize the material under consideration.
 Constants $y_0$ and $v_0$ are initial displacement and velocity.
 Operators $D^\gamma, \gamma \in [0,2]$ are operators of usual (for $\gamma$ integer)
 or fractional differentiation.

 System  $(CP)$ models various physical processes.
 For example, $(CP)_1$ represents an
 equation of motion of a physical pendulum (in this case $f(t,y(t)) = \sin y
 (t)$),
 moving in a dissipative media of viscoelastic type;
 $(CP)_2$ represents  a constitutive equation of a material
 where energy is dissipated. Since the upper bound in integrals
 in $(CP)_2$ is two, both viscoelastic and viscoinertial effects
 are  included.
 As experiments show, the presence of integral on the left
 hand side of $(CP)_2$ indicates that dissipation
 properties depend on the order of the derivative. The integral on the right
 hand side of $(CP)_2$ is a consequence of the well known principle of
 equipresence (cf. \cite{TruNoll:65}).

 Fractional operators $D^{\ga}$ and $I^{\ga}$, $\ga\in\R$,
 (see (\ref{defFD}) and (\ref{defFI}), respectively)
 are widely used in constitutive  equations
 for modeling viscoelastic bodies.
 A typical model is the generalized Zener model, where relations
  between stress and strain
 involve linear fractional differential operators.
 We refer to \cite{Atanackovic:02}, \cite{CapMai:71} and \cite{BagTorv:84} for details.
 A viscoelastic model of wave propagation
 with distributed order derivatives was analyzed in
 \cite{Caputo:67}.
 In \cite{Atanackovic:03} a distributed derivative model of viscoelastic body was proposed,
 and the restrictions which follow from the Second Law of Thermodynamics
 are stated.
 Distributed order derivatives appear in many other branches of
 physics. For example, Caputo in  \cite{Caputo:01} introduced distributed order
 fractional derivative to model dielectric induction and diffusion.
 In the series of papers \cite{CheGorSoc:02,CheGorSocGon:03,MP:07,MPG:06} the
 fractional diffusion equation, which was obtained from standard
 diffusion equations by replacing the first-order time derivative
 with the fractional derivative of order $\be\in (0,1)$ was analyzed.
 In \cite{BagTorv:00} and \cite{LorHart:03}  distributed derivative models
 were used in the context of
 system identification theory.

 Nahu\v sev, \cite{Nakhushev:03}, considers  distributed order
 derivatives of the form $\int_c^d D^{\al}yd\al$, $c,d \in\R$.
 He used the Laplace transform method to show the existence of a
 solution to equations of type $\int_c^0 D^{\al}yd\al=g$, $c<0$ and
 $\int_0^d D^{\al}yd\al=g$, $d>0$ in $L^1[a,b]$, with $g\in L^1 $.
  Recently, Kochubei \cite{Kochubei:07}
  gave the necessary condition on continuous function $\phi$
  in order that distributed order fractional differential equation of the form
  $\int_0^d \phi(\ga) D^{\al}y d\al=g$ has a solution.
 Some examples of mechanical models with distributed order
 fractional derivative were also studied in our previous works
 \cite{TAMBSP:05}, \cite{SPLJOTA:06} and \cite{TASP:04}.
 Those are all particular cases of system $(CP)$.
% For example, in \cite{SPLJOTA:06} we have considered $(CP)$ with $\phi_2(\ga)=\del(\ga)$.

 In order to solve system $(CP)$ we will use the Laplace transform of tempered
 distributions supported by the half line and
 the Schauder fixed point theorem which leads to a
 solution in an interval $[0,\delta],$ for some $\delta > 0.$
 Also, we will show existence and uniqueness for  %, for given $g\in\spp$,
 a distributed order fractional differential equation
 $\int_{\supp\phi}\phi(\ga)D^{\ga} z=g$  in $\spp$
 if $\phi$ is of the form  $\phi = \sum_{i=0}^k a_i\del(\cdot-\ga_i)$, $a_i\in \R,a_i\neq 0
 $, $\ga_i\in \R$, $i=0,1,...,k$
 so that this equation has a form of %(\ref{DOFDE}) reduces to
 a linear fractional differential equation
 $\sum_{i=0}^k a_iD^{\ga_i}z=g$. We refer to  \cite{KST:06} for solutions of
 linear fractional differential equation
 in spaces of continuous or locally integrable functions.
 Also we refer to \cite{KempSchBey:02}
 for solutions within spaces of distributions defined on the hole
 line obtained by the use of the Fourier transform method.

 In Section \ref{SC0}, we will define different types of solutions to
 $(CP)$ problem: {\it{classical, mild}} and {\it{ non-impact}} ones.
 Further on, in Section \ref{SecLFDE}, we state and prove
 the existence and uniqueness of a solution  to a general linear fractional
 differential equation in the frame of $\spp$, assuming that the corresponding
 holomorphic function $\sum_{i=0}^k a_i s^{\ga_i}$, $0\leq\ga_i < 2$, $\Re s>0$,
 does not have zeros. %This is condition $(A_0)$.
 The existence of solutions to system $(CP)$ is proved in Section \ref{SC1}.
 In Theorem \ref{Existence}, we give
 sufficient conditions for the existence of
 {\it {mild}}, {\it{non-impact}} and {\it{classical}} solution to $(CP)$,
 in the case when $\phi_1$ is continuous in $[c,d]$ for some $0\leq c < d<2$
 and $\phi _{2}(\cdot)=\sum_{i=0}^{k}a_{i}\delta (\cdot-\ga_{i})$, $0\leq\ga_i<2$.
 Further on, we analyze cases where form of $\phi_2$ remains the
 same, while $\phi_1$ changes.
 In  Section \ref{SC2}, we perform  our analysis for the case
 $\phi _{1}(\ga)=b^{\ga}$ and $\phi _{2}(\ga)=a^{\ga}$, $b>a$,
 where $b>a$ is a restriction which follows from the Second Law of Thermodynamics.
 %%%%%%%%%%%%%%%%%%%%%%%%%%%%%%%%%%%%%%%%%%%%%%%
 %\hspace{0.4cm}

 \subsection{Notation and notions}
 Let $I\subseteq \R$ be an interval of $\R$. Then, we denote
  the space of locally integrable functions in $I$ by
 $L^1_{\text{loc}}(I)$,
 the space of $k$-times continuously differentiable  functions $y$  by
 $C^k(I)$, and  the space of  functions $y\in C^{k-1}(I)$, such that $y^{(k)}\in
 L^1_{\text{loc}}(I)$ by $AC^k(I)$.

 The space of functions $y\in L^1_{\text{loc}}(\R)$,
 with the property $y(t)=0$, $t<0$  we will denote  by $\Lloc(\R)$.

 Further,
 $\mathscr{D}(\R)$ and $\mathscr{S}(\mathbb{R})$
 are the spaces of compactly supported smooth functions and
 of rapidly decreasing smooth functions in $\mathbb{R}$,
 respectively. Their duals
 $\mathscr{D}'(\R)$ and $\mathscr{S}^{\prime}(\mathbb{R})$ are
 the spaces of Schwartz's and tempered distributions, respectively.
 $\dpp (\R)$ and $\spp(\R)$ denote their
 subspaces consisting of distributions supported by
 $[0,\infty )$. Henceforth, we drop $\R$ in the notation
 of $\dpp (\R)$ and $\spp(\R)$.
 We will also use compactly supported distributions, elements of
 $\E'(\Omega)$, where $\Omega$ is an open interval of $\R$.
 Any element of $\E'(\Omega)$ can be extended to $\R$
 as an element of $\mathscr{D}'(\R)$ in the usual way.
 Let $J\subseteq I$ be a subinterval of $I$ and let $y$ be a function or
 a distribution defined in $I$. We will use the same
 notation $y$ for its restriction in $J$. Thus $y|_{J}$ will be denoted
 by $y$ again.
 In this sense, we will write $f\in \spp \cap C^p([0,a])$,
 which means that $f\in\spp$ and $f|_{[0,a]}\in C^p([0,a])$.

 The Riemann-Liouville operator of fractional differentiation or
 integration is denoted by $D^{\gamma }$, $\gamma \in \mathbb{R}$
 defined as follows.
 Let $y\in \Lloc(\R)$.
 The Riemann-Liouville fractional
 integral of order $\ga > 0$ is defined by
 \begin{equation}\label{defFI}
 { I^{\ga}y(t): = \frac{1}{\Gamma (\ga)}}\int_{0}^{t}(t-\tau )^{\ga-1}y(\tau )d\tau,\quad t >
 0,
 \end{equation}
 where $\Gamma $ is the Euler gamma function.
 For $\ga=0$ one defines $I^{0}y := y$.
 If $y\in \Lloc(\R)$, it is shown (cf. {\cite[Theorem 2.6]{SKM:93}})
 that
 $\lim_{\ga\to 0}I^{\ga}y(t)=y(t)$ almost everywhere in $\R$.

 %Let $\overline{\R}_+ = \{x\in\R;x\geq 0\}$.
 Let $y\in \Lloc(\R)$ and suppose that for every $a>0$, $y\in AC^k([0,a])$.
 The Riemann-Liouville fractional
 derivative of order
 $0 \leq \ga < k$ for some $k\in\N$,
 is defined by
 \begin{equation}\label{defFD}
  D^{\ga}y(t) := \frac{d^{k}}{dt^{k}}I^{k-\ga}y(t),\quad t > 0.
 \end{equation}
 If $\ga\in\N,$ then $D^{\ga}=\frac{d^{\ga}}{dt^{\ga}}$.

 Riemann-Liouville's fractional operators in the setting of distributions
 are defined (e.g. \cite{Vladimirov:79}) by introducing the family
 $f_{\al}\in\dpp$, $\al\in\R$
 \begin{equation*}%\label{falfa}
 f_{\alpha }(t) := \left\{
 \begin{array}{cc}
  H(t)\frac{t^{\alpha -1}}{\Gamma (\alpha )}, & t\in\R,\,\alpha > 0, \\
 \frac{d^{N}}{dt^{N}}f_{\alpha + N }(t), & \quad \alpha \leq 0,\,\alpha + N > 0 ,\, N \in
 \N,
 \
 \end{array}
 \right.
 \end{equation*}
 where $H$ is Heaviside's function. Then, %$\al\mapsto f_{\al}$ is weakly holomorphic
 $f_{\al}*f_{\be}=f_{\al+\be}$, for all  $\al,\be\in\R$ and
 $f_{\al}*$ is the convolution operator in $\dpp$. Also, $f_{\al}: \spp\to\spp$.
 The operator $f_{\alpha }\ast $ in $\dpp$ is the
 operator of fractional differentiation for $\alpha < 0$ and of
 fractional integration for $\alpha > 0$. It coincides with
 the operator of derivation for $-\al\in\N$ and
 integration for $\alpha \in \N_0 = \N \cup \{0\}$.
 In particular, $f_0=\delta$.

 If $\alpha > 0$ and $y\in \Lloc (\R)$ then $I^{\alpha }y = f_{\alpha }\ast y$.
 If $0 \leq \al < k $, $k\in\N$ and for all $a>0$,
 $y\in \Lloc (\R) \cap AC^k([0,a])$,  then $D^{\alpha}y = f_{-\alpha }\ast y$.
 We will use notion $D^{\alpha}y$ also when $y$ is distribution and
 will denote $D^{-\ga}y=I^{\ga}y$, $\ga > 0$. Note that
 $D^{\ga}I^{\ga}y = y$ for $y\in\Lloc(\R)$ and
 $I^{\ga}D^{\ga}y = y$, $\ga > 0$ if $y$ is considered as an element in $\dpp$.\\

 Let $0\leq c< d $, $\phi\in C([c,d])$
 and $y\in \Lloc (\R) \cap AC^2([0,a])$.
 Then the distributed order fractional derivative in $[c,d]$
 is given by
 \begin{equation}\label{dofd}
 \int_c^d \phi(\ga)D^{\ga}y(t)d\ga.
 \end{equation}
 If  $\phi:[c,d]\to \R_+$ be a piecewise continuous
 bounded function and $y\in\spp$, then as in \cite{TASP:04},
 one can  consider (\ref{dofd}) as a Bochner type of integral.
 We refer to \cite{ATOLJSP:08} for the properties of  (\ref{dofd}) within $\spp$.
 Recall the definition from \cite{ATOLJSP:08}:
 \begin{definition}%[\cite{ATOLJSP:08}]
 Let $\phi\in\E'(\R)$ and $y\in\spp$.
 Then $\underset{\supp \phi}{\int}\phi(\gamma)D^{\gamma}y\, d\gamma$
 is defined as an element of $\spp$ by
 \begin{equation*} %\label{DefInt}
 \langle \underset{\supp \phi}{\int} \phi(\gamma)D^{\gamma}y(t) \,d\gamma ,
 \varphi(t) \rangle: =
 \langle \phi(\gamma) , \langle D^{\gamma}y(t) , \varphi(t) \rangle \rangle,
 \quad \varphi\in\S (\mathbb{R} ).
 \end{equation*}
 \end{definition}
 \noindent When $\supp\phi\subset [a,b]$ we write
 $\int_a^b \phi(\gamma)D^{\gamma}y \,d\gamma$ instead of
 $\underset{\supp \phi}{\int}\phi(\gamma)D^{\gamma}y\, d\gamma$.
 It is shown in \cite{ATOLJSP:08} that  $\ga\mapsto\dis{D^{\ga}y}{\phi}:\R\to\R$ is smooth
 and that the mapping $(\alpha,y) \mapsto D^{\alpha}y$ is
 continuous from $\mathbb{R} \times\mathscr{S}^{\prime}_+$
 to $ \mathscr{S} ^{\prime}_+$. %, so (\ref{DefInt}) is well defined.
 Also,
 $
 y \mapsto \underset{\supp \phi}{\int}\phi(\gamma)D^{\gamma}y \,d\gamma
 $
 is a linear and continuous mapping from $\spp $ to  $\spp$.
 For
 \begin{align}\label{fisumdelta}
 \phi _{2}(\cdot)=\sum_{i=0}^{k}a_{i}\del(\cdot -\ga_{i}),\quad
 \ga_{i}\in \R\,i\in\{0,1,...,k\}
 \end{align}
 distributed order fractional derivatives becomes a sum of finite
 number derivatives of fractional order
 $\sum_{i=0}^{k}a_{i}D^{\ga_i}y$.

  Distributed order fractional differential equation is given by
  \begin{equation}\label{dofde}
  \underset{\supp \phi}{\int}\phi(\gamma)D^{\gamma}y \,d\gamma = g,
  \end{equation}
  while for $\phi$ of the form (\ref{fisumdelta}), it becomes linear
  fractional differential equation. % which will be studied in Section \ref{SecLFDE}.

 To deal with fractional differential equations, it is  convenient to use
 the Laplace transform.
 Recall, if $y\in \mathscr{S}'_+$, then its Laplace
 transform is defined by
 \begin{equation*}%\label{laplacetransform}
 \widehat{y}(s)=\mathcal{L}y(s)=\langle y(t),\varphi (t)e^{-st}\rangle ,\quad
 \Re {s}>0,
 \end{equation*}%
 where $\varphi \in C^{\infty }$, $\varphi =1$ in $(-a,\infty )$ and $\varphi =0$
 in $(-\infty ,-2a)$, $a>0$. Note that $\mathcal{L}y$ is an analytic function for
 $\Re s>0$, and that the definition of $\LT y$
 does not depend on a chosen function $\varphi $ with given
 properties. We will use the identities given in the next
 proposition (see \cite{ATOLJSP:08}).
 \begin{proposition}\label{prop1}
 Let $\phi\in\E'(\R)$ and $y\in\spp$. Then:\\ [0.1cm]
 a) $\LT({f_{\alpha}\ast y})(s) = \frac{1}{s^{\alpha }}\hat{y}(s),\, \Re s
 >0,\al\in\R.$\\ [0.2cm]
 b) ${\mathcal{L}}(\underset{\supp \phi}{\int} \phi(\gamma)D^{\gamma}y\,d\gamma)(s)
    = \hat{y}(s)\langle \phi(\gamma) , s^{\gamma} \rangle, \, \Re s> 0.$\\ [0.1cm]
 c) If $\phi $ is a continuous function in $[c,d]$ and $\phi(\ga)=0$,
 $\ga\notin [c,d]$, then
 $${\mathcal{L}}(\int_c^d \phi(\gamma)D^{\gamma}y\,d\gamma)(s)
   = \hat{y}(s)\int_c^d \phi(\gamma)s^{\gamma} \,d\gamma, \quad \Re s> 0.$$
 \end{proposition}

 \section{Definitions of solutions to $(CP)$}\label{SC0}

  Problem $(CP)$ will be analyzed in the context of the following
  definitions.
 \begin{definition} [Classical solution]\label{def1}
  Let $f: [0,\infty)\times\R\to\R $ be  continuous and $\phi_1,\phi_2\in\E'(\R)$,
  $\supp\phi_i\subset [0,2]$, $i=0,1$.
  Let $y_0,v_0\in\R$,
  $\delta>0$, $y\in \spp \cap C^2([0,\delta])$,
  $z\in \spp\cap C([0,\delta])$.
  % and $D^{\ga}z\in L^1([0,\del]) $,  $\ga\in \supp\phi_2$.

  We say that $(y,z)$ is a
  {\bf{classical solution to $(CP)$ in  $[0,\delta]$}}
  if  $y$ and $z$ satisfy $(\widetilde{CP})$, where
  \begin{equation*}
 (\widetilde{CP})\quad\quad
 \left\{ \begin{array}{c}
  D^{2}y(t)+z(t)=f(t,y(t)),\quad t\in [0,\del], \\
   \\
  \int_0^2 \phi_1(\gamma) D^{\gamma}y d\gamma =
  \int_0^2 \phi_2(\ga) D^{\gamma}zd\gamma,\quad \text{ in } \spp, \\
    \\
  y(0)=y_{0},\qquad y^{(1)}(0)=v_{0}.
  \end{array}
  \right.
  \end{equation*}
 \end{definition}
  Let $\phi_1,\phi_2\in\E'(\R)$, assume that $y\in\spp$ and let $z\in\spp$ be a solution to (\ref{dofde})
  with $\phi=\phi_2$ and $g=\int_{0}^2 \phi_1(\ga)D^{\ga}y d\ga$. Then
  $z$ is a solution to $(CP)_2$ in $\spp$.
  This is provided by
  \begin{equation}\label{inverz}
  \ILT\left(\frac{1}{\dis{\phi _{2}(\gamma
  )}{s^{\gamma}}}\right)\in\spp,
  \end{equation}
  and  then
  \begin{equation}\label{zet}
  z = \ILT\left(\frac{1}{\dis{\phi _{2}(\gamma )}{s^{\gamma}}}\right)
  \ast \int_{0}^{2}\phi _{1}(\gamma )D^{\gamma}y(\cdot)d\gamma, \quad \text{in }
  \spp.
  \end{equation}
  Equation $(CP)_1$  with initial conditions $(CP)_3$, in the setting of
  tempered distributions, reads
  \begin{equation}\label{ODEandIC}
  f_{-2}*y + z = f(\cdot, y) + y_0\del'+ v_0\del.
  \end{equation}
  We substitute (\ref{zet}) in (\ref{ODEandIC}) and apply $D^{-2}=f_2\ast$.
  It follows that
  $$
    y = - \ILT\left(\frac{1}{\dis{\phi _{2}(\gamma )}{s^{\gamma}}}\right)
    \ast \int_{0}^{2}\phi _{1}(\gamma )D^{\ga-2}y d\ga + D^{-2}f(\cdot,y))+ v_{0}f_1 + y_{0}H.
  $$
   \begin{definition}[Mild solution]\label{def3}
   Let $f: [0,\infty)\times\R\to\R $ be a locally integrable
   function and $\phi_1,\phi_2\in \E'(\R)$, $\supp\phi_i\subset [0,2]$, $i=0,1$,
   where $\phi_2$  is such that (\ref{inverz}) holds.
   Let $y_0,v_0\in\R$ and $\del>0$.
   The pair  $(y,z)$ is called a {\bf{mild solution to $(CP)$ in
   $[0,\del]$}} if $y\in\spp\cap C([0,\del])$, such that
   \begin{align}{\label{DDEC1}}
    y(t)=& - \ILT\left(\frac{1}{\dis{\phi _{2}(\gamma )}{s^{\gamma}}}\right)(t)
          * \int_0^2 \phi_1(\ga)D^{\ga-2}y(t)d\ga\nonumber \\
    &+ D^{-2}f(t,y(t))+ v_0t+ y_0,\quad t\in [0,\delta],
   \end{align}
   and $z$ is given by (\ref{zet}).
  \end{definition}
  \begin{remark}\rm
    a) A mild solution to $(CP)$ in $[0,\del]$ is a weak solution in $(0,\del)$.
    Namely,  (\ref{DDEC1}) and (\ref{zet})
    imply that for every
    $\theta \in C_0^{\infty}((0,\delta))$
    \begin{align*}
    \dis{D^2y}{\theta}& = \dis{y}{D^2\theta}
        = \dis{D^{-2}z(t) +  D^{-2}f(t,y(t))+ v_0t+ y_0}{D^2\theta(t)}\\
              & = \dis{z}{\theta} + \dis{f(t,y(t))}{\theta(t)}.
    \end{align*}
  b) A classical solution is also a mild solution.
  If $f: [0,\infty)\times\R\to\R$  is continuous and a mild solution
  exists, then considering that mild solution satisfies $(\widetilde{CP})_2$,
  in order for it to be a classical solution,
  one needs additionally to assume that
  $y\in C^2([0,\del])$ and $z\in C([0,\del])$.
  \end{remark}

  In mechanics,
  an impact solution is an absolutely continuous function
  whose first derivative has jumps.
  Impact solutions
  are included in the definition of the mild solution.
  %In this paper  a special attention is related to the observable quantity, this is strain $y$.
  In order to make an additional classification of solutions related
  to the observable quantity $y$, we
  introduce a class of {\it{non-impact solutions}}. This class will
  correspond to solutions in which the first derivative $y'$
  (velocity) is continuous.
  \begin{definition}[Non-impact solution]\label{Definicija_NonImpact}
  Let $f: [0,\infty)\times\R\to\R $ be a locally integrable function,
   $\phi_1,\phi_2\in \E'(\R)$, $\supp\phi_i\subset [0,2]$, $i=0,1$.
   Let $y_0,v_0\in\R$ and $\del>0$.
  The pair $(y,z)$ is called a {\bf{non-impact solution to $(CP)$ in
  $[0,\del]$}} if $y\in\spp\cap AC^2([0,\del])$ and $z\in\spp\cap
  L^1([0,\delta])$,
  such that
  $$
  D^{2}y(t)+z(t)=f(t,y(t)),\quad \text{almost everywhere }\quad t\in [0,\del],
  $$
  $$
  \int_0^2\phi_1(\gamma) D^{\gamma}y(t) d\gamma =
  \int_0^2 \phi_2(\ga) D^{\gamma}z(t)d\ga,\quad \text{in }\spp
  $$
  and
  $$
  y(0)=y_0,\quad Dy(0)=v_0.
  $$
  \end{definition}

 It is clear that condition (\ref{inverz}) plays an important role in definitions of the
 solutions, and that satisfaction of (\ref{inverz}) depends on
 $\phi_2$. However, in the next section, we will see
 that when $\phi_2$ is of the form (\ref{fisumdelta})
 with additional condition $(A_0)$ satisfied, then (\ref{inverz}) is fulfilled.
 As we will see in Section \ref{SC1}, (\ref{inverz}) and certain
 conditions on $f$  assure the existence of mild, non-impact, and classical
 solutions.

 \section{Linear fractional differential equation in $\spp$}\label{SecLFDE}

 In this section, we
 prove the existence and the uniqueness of a  fundamental solution  to a general
 linear fractional differential equation in $\spp$.

  Let $g\in\spp$. Let $a_i\neq 0$ and $\ga_i$ be  arbitrary real
  numbers,
  such that $\ga_i > \ga_{i+1}$, $i\in\{0,1...,k\}$,
  and consider
  \begin{equation}  \label{CEC11}
   \sum_{i=0}^k a_i D^{\gamma_i}z = g,\quad {\text {in } }\,\spp.
  \end{equation}

  %Our main result of this section is the next theorem.
 \begin{theorem}\label{T1}
  Equation (\ref{CEC11}) has a unique solution $z\in\spp$, if and only if
  \begin{align*}
   (A_0)\quad\quad \sum_{i=0}^k a_i s^{\gamma_i}\neq 0,\quad s\in\mathbb{C} _+ = \{ s\in\C\,;\, \Re s >0\}.
  \end{align*}
  \end{theorem}
  %%%%%%%%%%%%%%%%%%%%%%%%%%%%%%%%%%%%%%%%%
  \pr %{\it Proof of Theorem \ref{T1}}.
 %%%%%%%%%%%%%%%%%%%%%%%%%%%%%%%%%%%%%%%%%
  First note
  \begin{equation*}%\label{conformLFDE}
   \sum_{i=0}^k a_i D^{\ga_i}z = \sum_{i=0}^k a_i (f_{-\ga_i}*z) =
   (\sum_{i=0}^k a_if_{-\ga_i})*z,
  \end{equation*}
  i.e. that (\ref{CEC11}) is equation of a convolution type.
  For given $\varphi, h\in\spp$, the equation $\varphi * z = h$ is uniquely
  solvable    if and only if there exists $a,b\in%
 \mathbb{R} $ and $C>0$, such that
 \begin{equation*}%\label{ISO}
  \frac{1}{|\widehat{\varphi}(s)|} \leq C\frac{(1+|s|)^a}{|\Re s|^b},\quad
  s \in\C_+ ,
 \end{equation*}
 where $\widehat{\varphi}$ denotes Laplace transform of $\varphi$,
 (see \cite[Chap.7, p.50]{VlaDroZav}).

 Let us denote $$F(s)= \wh{\vphi}(s) = \sum_{i=0}^k a_i s^{\ga_i},\,s\in\C_+.$$
 We will show that $F(s)\neq 0$, $ s\in \C_+$ (which is $(A_0)$)
 implies that there exist $a\in\R$ and $C>0$, such that
 \begin{equation}\label{broj}
 \frac{1}{|F(s)|}\leq C \frac{(1+|s|)^a}{(\Re s) ^a},\, s\in\C_+.
 \end{equation}
 This will imply the theorem.

 To begin with, note that function $F$ is analytic in $\C \backslash (-\infty,0]$.
 Next observe that
 $|F(s)|$ behaves like $|a_0s^{\ga_0}|$ when $|s|$ is large, and like
 $|a_ks^{\ga_k}|$ when $|s|$ is small which implies that there are $r,R\in\R$,
 $0<r<R$, such that $F(s)$ has no zeros for $|s|<r$ and $|s|>R$.
 Further observe that $F(s)\neq 0$ for $\Re s>0$ implies that $F(s)\neq 0$ for $\Re
 s<0$ and therefore the  set of zeros
 of function $F$ is contained in $[-iR,-ir]\cup [ir,iR]$.
 Thus, it must be a finite set of points $iy_j$, $j\in \{1,2...l\}$,
 since zero set of an analytic function does not have limits in the
 domain of its analyticity, so in any compact region,
 the set of zeros can only be a finite set.

 Since $F$ is analytic, it follows that in  neighborhood
 of zeros $iy_j$,
 $|F(s)|$ behaves as $|s-iy_j|^{m_j}$, $m_j\in\N$, $j\in\{1,2,..,l\}$.
 Denote $D_j=\{s\in\C; |s-iy_j|<r_i\}$, for some $r_j>0$, $j\in
 \{1,2...l\}$, $D=\bigcup_{j=1}^{ l} D_j$ and $K=\{s\in\C; r\leq |s| \leq R\}$.
 Choose $a\geq m_j$, $j\in \{1,2...,l\}$
 and $a\geq \ga_0$. There exist $C$ such that
 $$
 \frac{1}{|F(s)|}\leq C \frac{1}{(\Re s) ^a},\, s\in
 (\C\backslash K) \cup D.
 $$
 In $K\backslash D$ function $F$ reaches its minimum, denoted by $m$, so
 with some new $C$
 $$
 \frac{1}{|F(s)|}\leq  \frac{1}{m}\leq C \frac{(1+|s|)^a}{(\Re s) ^a},\, s\in
 K.
 $$
 Putting all together we arrive to (\ref{broj}).
 \ep

  \begin{remark}\rm\label{MechRem1}
  Equation (\ref{CEC11}) represents a constitutive equation of a
  viscoelastic body which, coupled with equation
  $(CP)_1,$ represents a mechanical model $(CP)$.
  We explained in \cite{ATOLJSP:08} the mechanical aspect of
  condition ($A_0$) and here we repeat it briefly. Let $T>0$ and
  $y$ be of the form $y(t)=\sum_{i=0}^{k}a_{i}D^{\ga_{i}}z(t)$, $t\in
  [0,T]$. The dissipation work is given by
  $A_d=\int_{0}^{T}z(t)y^{(1)}(t)dt$ and the
  dissipation inequality reads $A_d\geq 0$ (see \cite{Christensen:82}).
  If $(A_0)$ does not hold, then the dissipation inequality  is violated.
  So, $(A_0)$ is a necessary condition if $(CP)$
  represents a mechanical model.
 \end{remark}

 %%%%%%%%%%%%%%%%%%%%%%%%%%%%%%%%%%%%%%%%%

 %Applying Laplace transform to (\ref{conformLFDE}) one can see that
 The solution to (\ref{CEC11}) is of the form
 $z=l*g$, where
 $
 l% = \ILT\left(\frac{1}{\sum_{i=0}^k a_i s^{\ga_i}}\right)
 $
 is a fundamental solution to (\ref{CEC11}), i.e. solution to
 $\sum_{i=0}^k a_i D^{\ga_i}y=\delta$, with
 the properties  given in next lemma. For the proof see \cite{ATOLJSP:08}.
 \begin{lemma}\label{mg}
  Let $\ga_i\in [0,2)$ and $\ga_i > \ga_{i+1} \geq 0$ for all $i\in\{0,1...,k\}$.
  Assume $(A_0)$. Let
  \begin{equation}\label{l}
   l(t)=\left\{ \begin{array}{cc}
            \ILT\left(\frac{1}{\sum_{i=0}^k a_i s^{\ga_i}}\right)(t),& t\geq 0 \\
            0 \quad,& t<0.
            \
        \end{array}
        \right.
  \end{equation}
 Then:\\
  %\begin{itemize}
 (i) $l$ is a locally integrable function in $\R$.\\
 (ii) $l$ is absolutely continuous in $\R$, if $\ga_0 - \ga_k > 1$.
 %  \end{itemize}
 \end{lemma}

 \section{Solvability of the system $(CP)$} \label{SC1}

 Through this section we consider
 cases when $\phi_2$ is of the following form
  \begin{align}\label{fidva}
 \phi _{2}(\cdot)=\sum_{i=0}^{k}a_{i}\del(\cdot -\ga_{i}),\quad
 \ga_{i}\in [0,2),\,\ga_{i}> \ga_{i+1},\,i\in\{0,1,...,k\}.
 \end{align}
 The constitutive equation $(CP)_2$ then reads
 \begin{equation}  \label{CEC1}
  \int_0^2 \phi_1(\gamma) D^{\gamma}y(t) d\gamma = \sum_{i=0}^k a_i
  D^{\gamma_i}z(t), \quad t>0.
 \end{equation}
  If condition $(A_0)$ of Theorem \ref{T1} is satisfied,
  and $l$ is defined as in Lemma \ref{mg} by (\ref{l}), then,
  if we suppose $y\in\spp$,
  Theorem \ref{T1} implies that (\ref{CEC1}) has a unique
  solution $z\in\spp$ given by
  \begin{equation*}%\label{z}
   z = l *   \int_0^2 \phi_1(\ga) D^{\ga}y d\ga.
  \end{equation*}
  Equation
  (\ref{DDEC1}) becomes
  \begin{align}\label{DDEC2}
    y(t)= - l * \int_0^2 \phi_1(\ga)D^{\ga-2}y(t)d\ga + D^{-2}f(t,y(t))+ v_0t+ y_0,\quad
    t>0,
   \end{align}
  so a mild solution to $(CP)$ will exist if (\ref{DDEC2}) has a solution $y\in
  C([0,\del])\cap \spp$, for some $\del>0$.
  To show such existence, we will define  a  mapping  $T$
  by the right hand side of (\ref{DDEC2}),
  and show existence of its fixed point in appropriate Banach space.
  The necessary conditions for $T$ to have a fix point are to
  follow.\\

  \noindent First, let us impose conditions on $f$:
 \begin{align*}
 (A_l)\quad & f(t,u),\,(t,u)\in [0,\infty )\times \R, \text{ is locally integrable and }\\
 & \text{ there exist }a>0, \, \alpha >0\text{ and }h\in L^{1}([0,a]),\,h \geq 0,
 \text{ such that } \\
 & |f(t,u)|\leq h(t)|u|^{\alpha },\, u\in \mathbb{R},t\in [0,a],
 \end{align*}
 and the stronger one
 \begin{align*}
 (A_l)'\quad & f(t,u),\,(t,u)\in [0,\infty )\times \R, \text{ is continuous }
  f(t,0)=0, t\in[0,\infty),\text{ and} \\
  & \text{ there exist }\al>0,\, a>0,\, h\in C([0,a])\text{ with } h\geq 0 \text{ and }h(0)=0,
  \text{ such that } \\
  & |f(t,u)-f(t,v)|\leq h(t)|u-v|,\,u,v\in \R,t\in [0,a].
 \end{align*}
 Further, we introduce  conditions on $\phi_1$ and $\phi_2$.
  \begin{align*}
  (\Phi_1) \quad
  &\phi_1 \text{  is continuous in  $[c,d]\subset [0,2)$, $\phi_1(\ga)=0$, $\ga\notin
  [c,d]$};\\
  &\phi _2 \text{ is of the form (\ref{fidva})}.
  \end{align*}
 The following theorem holds.
  \begin{theorem}\label{Existence}
  Let $y_0,v_0\in \R$, $r>\max\{|y_0|,|v_0|\}$.
  Let $(\Phi_1)$, $(A_0)$ and $(A_l)$ hold.
  There exists $\delta=\del(r) > 0$
  such that:
  \begin{itemize}
   \item[a)] $(CP)$ has a mild solution in $[0,\del]$.
   \item[b)] If  $y(0)=y_0=0$, then $(CP)$ has a mild solution $(y,z)$ in
   $[0,\del]$,
    such that $y\in C^1([0,1])$. If $\ga_0-\ga_k> 1$ then $(y,z)$ is a non-impact solution.
   \item [c)] If $(A_l)'$ holds and $y_0 = v_0 = 0$,
 then $(CP)$ has a classical solution in $[0,\del]$.
 \end{itemize}
 \end{theorem}
 \pr
  We will give the joint proof for all parts of the theorem, underlying
  the differences if we assume  $(A_l)$ or $(A_l)'$.

  Let $l$ be defined by (\ref{l})  as in Lemma \ref{mg}
  and  define operator $T$ by the right-hand side of (\ref{DDEC2}),
  i.e.
  $$
   Ty(t):=- l * \int_0^2 \phi_1(\ga)D^{\ga-2}y(t)d\ga + D^{-2}f(t,y(t))+ v_0t+
   y_0.
  $$
  We will consider $T$  acting on Banach spaces
   $C^i([0,a])$ endowed  with the norms
     $\|y\|_i := \sup_{t\in [0,a]}\{|y^{(j)}(t)|; 0\leq j\leq i\}$,
  $i=0,1,2$ and, in particular, Banach spaces
  $$C^i_{0}([0,a]):=\{y\in C^i([0,a]); y(0)= 0\} $$
  and
  $$
  C^i_{00}([0,a]) := \{y\in C^2([0,a]);y(0)= 0, y'(0)=0\}
  $$
  endowed with the same norm. We will show that
  $T$ is a compact operator in quoted spaces and that $T$ maps a closed ball into itself.
  Then we will apply the Schauder fixed point theorem and obtain a mild
  solution.\\
  %in one of the quoted spaces depending whether we have assumed $(A_l)$ or
  %$(A_l)'$. \\%(Of course if something holds with $(A_l)$ it holds with $(A_l)'$).\\

 (i) {\it $T$ maps $C^i([0,a])$ into $C^i([0,a])$ for $i=0,1,2$.}\\

  Let $y\in C([0,a])$, where $a$ comes from $(A_l)$ or $(A_l)'$.
  Define
  \begin{align}\label{J}
   Jy(t) &:= \int_c^d \phi_1(\ga) D^{\ga-2}y(t)d\ga \nonumber\\
   &=\int_{c}^{d}\frac{\phi_1(\ga)}{\Gamma(2-\gamma)}
   \int_{0}^{t}(t-x)^{1-\gamma}y(x)dxd\gamma,\quad t\in[0,a].
  \end{align}
  Since $y\in C([0,a])$ and
    $\supp\phi_1\subset [c,d]\subset [0,2)$,
  $Jy$ is continuous. % by its definition (\ref{J}).
 Further, denote
  \begin{equation*}%\label{R}
   Ry(t) := (l*Jy)(t)= \int_0^t l(x)Jy(t-x)dx,\quad t\in[0,a],
  \end{equation*}
  where $l$ is defined as in Lemma \ref{mg} and therefore is  locally
  integrable. It follows that $Ry$ is continuous. % by its definition (\ref{R}).
  Function  $Gy$ defined by
  \begin{equation*}%\label{G}
   Gy(t):= D^{-2}f(t,y(t))=\int_0^t (t-x) f(x,y(x))dx,\quad t\in[0,a]
  \end{equation*}
  is continuous by assumption $(A_l)$.
  Since
  \begin{equation}\label{opT}
   Ty(t):= Ry(t)+Gy(t)+ v_0 t + y_0,\quad t\in[0,a],
  \end{equation}
  it follows $T:C([0,a])\to C([0,a])$.\\

  In what follows we shell need well known
  fact stated in next lemma.
  \begin{lemma}\label{lemmaL}
  (a) Let $l\in C(\R)$ and $\supp{l}\subset [0,\infty)$.
  Then the convolution operator $l*$,
  defined by $y\mapsto l*y$, for $y\in\Lloc(\R)$, maps continuously
  $C^1([0,a])$ into $C_0^1([0,a])$
  and $C^2_0([0,a])$ into $C_{00}^2([0,a])$.

  (b) Let  $l\in\Lloc(\R)$. Then  $y\mapsto l*y$, $y\in\Lloc (\R)$
  is a continuous mapping
  from $C([0,a])$ into $C_0([0,a])$,
  from $C^1_{0}([0,a])$ into $C^1_{00}([0,a])$ and
  from $C^2_{00}([0,a])$ into $C^2_{00}([0,a])$.\\
  \end{lemma}

  Let  $y\in C_0^1([0,a])$. We will show that $Ty\in C_0^1([0,a])$.
  Partial integration in (\ref{J}) and $y(0)=0$ follows
  \begin{align}\label{prviizvodJ}
  \frac{d}{dt}Jy(t)& =\int_c^d\frac{\phi_1(\ga)}{\Ga(2-\ga)}
   \int_0^t(t-x)^{1-\ga}y^{(1)}(x)dxd\ga
    +  y(0)\int_c^d\frac{\phi_1(\ga)}{\Ga(2-\ga)}t^{1-\ga}d\ga,\nonumber \\
    & = Jy^{(1)}(t) +  y(0)\int_c^d\frac{\phi_1(\ga)}{\Ga(2-\ga)}t^{1-\ga}d\ga, \quad   \nonumber \\
    & = Jy^{(1)}(t)\quad t\in [0,a].
  \end{align}
  Therefore $\frac{d}{dt}Jy$ is continuous and $ J:C_0^1([0,a])\to
  C^1_0([0,a])$. This, and the fact that $l$ is locally integrable, by Lemma
  \ref{lemmaL}, implies that  $R$ maps $C_0^1([0,a])$ into $C_0^1([0,a])$ and
  \begin{equation*}%\label{prviizvodR}
   \frac{d}{dt}Ry(t) = l*\frac{d}{dt}Jy(t) =  l*Jy^{(1)}(t),\quad  t\in [0,a].
  \end{equation*}
  By $(A_l)$ it follows that
  \begin{equation*}%\label{prviizvodG}
  \frac{d}{dt}Gy(t)= \int_0^t f(x,y(x))dx, \quad t\in [0,a],
  \end{equation*}
  is continuous and $G:C_0^1([0,a])\to C_0^1([0,a])$.\\
  In total we have that
  \begin{equation}\label{prviizvodT}
  \frac{d}{dt}Ty(t)= l*Jy'(t) + \int_0^t f(x,y(x))dx + v_0,\quad
  t\in[0,a]
  \end{equation}
  is continuous, $Ty(0)=0$,  and  therefore $T:C_0^1([0,a])\to C_0^1([0,a])$.

  Finally, let $y\in C_{00}^2([0,a])$.
  We will show that $Ty\in C_{00}^2([0,a])$.
  Since $y(0)=0$ and $y'(0)=0$
  \begin{align*}%\label{drugiizvodJ}
  \frac{d^2}{dt^2}Jy(t) = Jy^{(2)}(t) &+
    y'(0)\int_c^d\frac{\phi_1(\ga)}{\Ga(2-\ga)}t^{1-\ga}d\ga\nonumber\\
    &+ y(0)\int_c^d\frac{\phi_1(\ga)(1-\ga)}{\Ga(2-\ga)}t^{-\ga}d\ga, \quad
    t\in[0,a]
  \end{align*}
  is continuous. Thus,
  \begin{equation}\label{IzvodiJ}
  \frac{d^i}{dt^i}Jy(t)=Jy^{(i)}(t), \quad t\in[0,a],
  \,i\in\{0,1,2\},
  \end{equation}
  $Jy(0)=Jy'(0)=0$, so  $J:C_{00}^2([0,a])\to C_{00}^2([0,a])$.
  This, with Lemma \ref{lemmaL}, %(with $Jy(0)=0$ and  $Jy'(0)=0$),
  and the fact that $l$ is locally integrable, implies that
  $R$ maps $C_{00}^2([0,a])$ into $C_{00}^2([0,a])$ and
  \begin{equation*}%\label{drugiizvodR}
   \frac{d^i}{dt^i}Ry(t) = l*\frac{d^i}{dt^i}Jy(t)=
    l*Jy^{(i)}(t),\quad t\in[0,a],\quad i\in\{0,1,2\}.
  \end{equation*}
  Assume now that $(A_l)'$ holds. Then
  \begin{equation*}%\label{drugiizvodG}
   \frac{d^2}{dt^2}Gy(t)= f(t,y(t)),\quad t\in [0,a]
  \end{equation*}
  is continuous, $Gy(0)=\frac{d}{dt}Gy(0)=0$ and therefore $G:C_{00}^2([0,a])\to C_{00}^2([0,a])$.
  Thus,
  $T: C_{00}^2([0,a])\to C_{00}^2([0,a])$ and
  \begin{equation}\label{drugiizvodT}
   \frac{d^2}{dt^2}Ty(t)= l*Jy''(t) +  f(t,y(t)),\quad t\in[0,a].
  \end{equation}\\

  (ii) {\it Compactness of $T$.}\\

  Let us show that $J$ is a compact operator in $C([0,a])$.
  Recall,
  \begin{align*}
  Jy(t)& = \lim_{N\rightarrow\infty}\sum_{n=1}^{N}\phi_1(\gamma_{n})
  D^{\gamma_{n}-2}y(t)\Delta\gamma_{n}\\
       & = \lim_{N\rightarrow\infty}\sum_{n=1}^{N}\frac{\phi_1(\gamma _{n})}
       {\Gamma(2-\gamma_{n})}\int_{0}^{t}\frac{y(x)}{(t-x)^{\gamma_{n}-1}}dx\Delta\gamma_{n},
       \quad t\in [0,a],
 \end{align*}
 where $\gamma_{n}$ are points of the interval $[c,d]$ in the
 definition of the Riemann sum defined for the integral.
 Since $\ga_n-2\leq d-2<0$,  for all
 $n\in\N$, operators $ D^{\gamma_{n}-2}: C^i([0,a])\to C^i([0,a])$,
 $i=0,1,2$, $n\in\{1,...,N\}$ are operators of fractional integration and therefore compact.
 $J$ is a compact operator as a limit of compact operators
 in the operator norm topology.
 By the same argument, by (\ref{IzvodiJ}),
 $J: C_0^1\to C_0^1$ and $J: C_{00}^2\to C_{00}^2$ are compact.

 By Lemma \ref{lemmaL}, $l*$ is a continuous mapping: $C([0,a])\to C([0,a])$,
 $C_{0}^1([0,a])\to C_{0}^1([0,a])$ and  $C_{00}^2([0,a])\to
 C_{00}^2([0,a])$,
 respectively,
 and therefore $R$ is a compact operator
 as a composition of a compact and a continuous operator
 on $C([0,a])$, $C_{0}^1([0,a])$ and  $C_{00}^2([0,a])$, respectively.

 Under the assumption $(A_l)$, $G:C^1_0([0,a])\to C^1_0([0,a])$ is compact.
 Let $ M>0$ and $B_M$ be the ball in $C^1_0([0,a])$ ($B_M := \{y\in C^1_0([0,a]) ;  \|y\|_1\leq M \}$).
 The set $ G[B_M]:=\{G\theta ; \theta\in B_M\}$ is uniformly bounded.
 Let us prove this. Condition $(A_l)$ implies
 $$
 |Gy(t)|\leq \int_0^t |t-x||f(x,y(x))|dx\leq a M^{\al}C, \quad t\in [0,a],\quad  y\in B_M,
 $$
 and
 $$
 |\frac{d}{dt}Gy(t)| \leq \int_0^t |f(x,y(x))|dx \leq M^{\al}C, \quad t\in [0,a],\quad y\in B_M,
 $$
 where $\quad C=\max_{t\in [0,a]}\int_0^t h(x)dx$.
 It is also an equicontinuous family since
 \begin{align*}
 |Gy(t_1)&-Gy(t_2)|\leq\\
  &\leq \int_0^{t_1}|(t_1-x)-(t_2-x)||f(x,y(x))|dx
  + \int_{t_1}^{t_2}|t_2-x||f(x,y(x))|dx\\
  & \leq M^{\al}C(t_1-t_2) + \sup_{x\in [t_1,t_2]}|t_2-x|M^{\al}C \\
  &\leq (t_1-t_2)M^{\al}C(1+a)\leq C'(t_1-t_2),\quad t_1,t_2\in[0,a]
 \end{align*}
 and
 \begin{align*}%\label{equicontDG}
 |\frac{d}{dt}Gy(t_1)-\frac{d}{dt}Gy(t_2)|& \leq
 \int_{t_1}^{t_2}|f(x,y(x))|dx\nonumber\\
    &\leq M^{\al}\int_{t_1}^{t_2}|h(x)|dx ,\quad t_1,t_2\in[0,a].
 \end{align*}
 Since the function $t\mapsto\int_{0}^{t}|h(x)|dx$,
 $t\in[0,a]$ is absolutely continuous in $[a,b]$, it follows that
 the set $G[B_M]$ is equicontinuous.
 Therefore,  Arzela-Ascoli theorem implies that $G[B_M]$ is relatively compact and
 $G$ is, therefore, a compact operator in $C^1_0([0,a])$.

 If we use the stronger condition
 $(A_l)'$, then $G$ is compact in $ C^2_{00}([0,a])$.
 We will show  that $G$ maps a bounded  sequence $\{y_n\}_{n\in\N}$ in $C^2_{00}([0,a])$ to
 a relatively compact set $\{Gy_n\}_{n\in\N}$ in $C^2_{00}([0,a])$.
 If $\{y_n\}_{n\in\N}$ is bounded in $C^2_{00}([0,a]) $, then $\{y_n\}_{n\in\N}$ and
 $\{y_n'\}_{n\in\N}$ are equicontinuous and Arzela-Ascoli theorem assures that we have
 convergent subsequences, again
 denoted by $\{y_n\}_{n\in\N}$ and $\{y_n'\}_{n\in\N}$.
 Denote $ c_n = Gy_n $, $n\in\N$. Then, for $ i\in \{0,1\}$,
  \begin{align}\label{G01}
   |\frac{d^i}{dt^i}c_n(t)-\frac{d^i}{dt^i}c_m(t)|
   &\leq \int_0^t|t-x|^{1-i}|f(t,y_n(t))-f(t,y_m(t))|\nonumber \\
   &\leq D_i|y_n(t)-y_m(t)|^{\al}, \quad t\in [0,a],
  \end{align}
  where
  $D_i = \max_{t\in [0,a]} \int_0^t |t-x|^{1-i}|h(x)|dx,\,
  i\in\{0,1\},$ and
  \begin{align}\label{G2}
    |\frac{d^2}{dt^2}c_n(t)-\frac{d^2}{dt^2}c_m(t)|
    &\leq |f(t,y_n(t))-f(t,y_m(t))|\nonumber \\
    &\leq D_2|y_n(t)-y_m(t)|^{\al}, \quad t\in [0,a],
  \end{align}
 where $D_2 = \max_{t\in [0,a]}|h(x)|dx.$
 Since $\{y_n\}_{n\in\N}$
 converges, (\ref{G01})
 and (\ref{G2}) imply that the same is true
 for $\{c_n\}_{n\in\N}$ in $C^2([0,a])$.

 As a sum of compact operators $R$, $G$
 and the mapping: $y\mapsto v_0t + y_0$,
 $T$ is  compact and $T: C[0,a]\to C[0,a]$.
 Moreover,  $T: C^1_0([0,a])\to  C^1_0([0,a])$ is compact
 if $(A_l)$ holds, while
 $T: C^2_{00}([0,a])\to C^2_{00}([0,a])$ is  compact if $(A_l)'$ holds.\\

 (iii) {\it Determination of $\del$.}\\

  We will show that there exists $\del>0$ depending on $r$ such that
  $T$ maps  $B_r \subset C([0,\delta])$ ,
  $B_r\subset C^1_{0}([0,\delta])$ and $B_r\subset C^2_{00}([0,\delta])$ into
  itself.% ($B_r$ is the closed ball with radius $r$ and center at $0.$)

  From (\ref{J}) and  (\ref{IzvodiJ}) we  derive
  \begin{equation}\label{ocenaJ}
  |\frac{d^i}{dt^i}Jy(t)|\leq M_t \sup_{x\in [0,t]}|y^{(i)}(x)|, \quad i\in\{0,1,2\},\quad
  \end{equation}
  where
  \begin{equation*}
  M_t = \int_c^d \frac{|\phi_1(\ga)|}{\Ga(3-\ga)}t^{2-\ga}d\ga\quad t\in [0,\del].
  \end{equation*}

  Now consider $B_r\subset C([0,\delta])$ and $B_r \subset C_0^1([0,\delta])$,
  respectively.
  Then (\ref{opT}), (\ref{prviizvodT}), (\ref{ocenaJ}) and
  $(A_l)$ imply that for $i = 0,1$,
  \begin{align}\label{ocenaT}
  |&\frac{d^i}{dt^i}Ty(t)|\leq |l*\frac{d^i}{dt^i}Jy(t)| +
        |\frac{d^i}{dt^i}Gy(t)| + |v_0|t^{1-i} + |y_0|(1-i)\nonumber\\
  &\leq \int_0^t|l(t-x)||Jy^{(i)}(x)|dx + \int_0^t |(t-x)^{1-i}||f(x,y(x))|dx + |v_0|t^{1-i} + |y_0|(1-i)\nonumber \\
  &\leq D_t M_t\sup|y^{(i)}(t)| + C_t\sup|y(t)|^{\al}+ |v_0|t^{1-i} + |y_0|(1-i), \quad
  t\in[0,\del],
  \end{align}
  where $$ D_t = \int_0^t|l(x)|dx,\quad C_t = \int_0^t h(x)(t-x)^{1-i}dx,\quad t\in [0,\del].$$
  Since $|y^{(i)}(t)|\leq r$, $i=0,1$ and $r > \max\{|v_0|,
  |y_0|\}$,
  it is possible to shrink $\del$ as much as
  we need to reach
  $$ |\frac{d^i}{dt^i}Ty(t)|\leq  r, \quad \text{ for } i=0,1. $$
  With this, we determine $\delta$  in cases a) and b).

 For the case c),  consider the ball $B_r \subset C^2_{00}([0,\delta])$.
 Then (\ref{drugiizvodT}), (\ref{ocenaJ}) and  $(A_l)'$ imply that the estimates for
 $\frac{d^i}{dt^i}Ty(t)$, $i=0,1$ are
 valid  and
 \begin{align*}
  |\frac{d^2}{dt^2}Ty(t)|&\leq |l*\frac{d^i}{dt^i}Jy(t)| + |f(t,y(t))|\\
  &\leq D_t M_t\sup|y^{(2)}(t)| + C_t|y(t)|^{\al}, \quad t\in [0,\del],
 \end{align*}
 with $C_t = \max_{t\in[0,\del]}|h(t)|$  and $D_t$ and $M_t$ as above.
 Since $|y^{(i)}(t)|\leq r,\, i=0,1,2$ and $ r>0$, it is
 possible to shrink $\del$ such that
  $$ |\frac{d^i}{dt^i}Ty(t)|\leq  r, \quad \text{ for } i=0,1,2, $$
 so it determines $\del$. \\

 (iv) {\it{Assertions a), b) and c).}}\\

  $a)$ Recall that assuming $(A_l)$ and fixing
  $r>\max\{|y_0|,|v_0|\}$,
  by (i)-(iv), we have obtained the
  existence of $\del=\del(r)>0$, such that
  $T$ is a compact mapping $C([0,\del])\to C([0,\del]) $
  and that $T$ maps $B_r\subset C([0,\delta])$ into  itself.
  Thus, according to the Schauder fixed point theorem,
  $T$ has a fixed point in $B_r\subset C([0,\del])$,
  i.e. (\ref{DDEC2}) has a solution
  $y\in B_r\subset C([0,\delta])$.
  %By (\ref{absT}) it follows that
  %$y\in AC([0,\delta])$.
   Set
  \begin{equation}\label{z1}
   z = l * \int_c^d\phi_1(\ga)D^{\ga}y d\ga,\quad\text{with }
   l=\ILT(\frac{1}{\sum_{i=0}^k a_is^{\ga_i}}).
  \end{equation}
  Then $(y,z)$ is a mild solution to $(CP)$.

  $b)$ Again, assuming $(A_l)$ and  fixing $r$,  by (i)-(iv),
  we find  $\del=\del(r) > 0$, such that
  $T: C^1_0([0,\del])\to C^1_0([0,\del])$ is compact
  and  $T$ maps the ball $B_r\subset C^1_0([0,\delta])$ into itself.
  Again, the use of the Schauder fixed point theorem assures that
  $T$ has a fixed point $y$ in $B_r\subset C_0^1([0,\del])$. If $z$ is  given by (\ref{z1})
  then $(y,z)$ is a mild solution to $(CP)$ in $[0,\del]$.

  If $\ga_0-\ga_k > 1$, then Lemma \ref{mg} implies that $l$ is absolutely
  continuous. Thus, $l*Jy\in AC^2([0,a])$. Also $Gy\in
  AC^2([0,\del])$. By (\ref{opT}) and the fact that $y=Ty$ we obtain $y\in AC^2([0,a])$.

  $c)$ With $(A_l)'$ and $r>0$, by (i)-(iv),
  we find  $\del=\del(r) > 0$ such that
  $T$ is  compact and maps $B_r\subset C^2_{00}([0,a])$ into itself.
  According to the Schauder theorem, (\ref{DDEC2}) has a solution $y$ in
  $C^2_{00}([0,\del])$. If  $z$ is given by (\ref{z1}) then $(y,z)$
  is a mild solution to $(CP)$ in $[0,\del]$.
  $D^{\ga}y$, $\ga<2$ is continuous, since $y\in C^2_{00}$, and
  $z\in C([0,\del])$ by (\ref{z1}). Applying $D^2$ to $y$ given by
  (\ref{DDEC2}) leads to conclusion that $(y,z)$ satisfies
  $(\widetilde{CP})_1$.\\
  \ep

 \noindent {\textbf{{Different conditions on $\phi_1$ and  $f$
 }}}

  Further, we assume restriction on  $\supp \phi_1$, which
  will imply better regularity properties
  of solutions. Let
  \begin{align*}
   (\Phi_2) \quad
   &\phi_1 \text{ is continuous in $[c,d]\subset [0,1)$, $\phi_1(\ga)=0$, $\ga\notin
   [c,d]$};\\
   &\phi _2 \text{ is of the form (\ref{fidva})}.
  \end{align*}
  Then, we have the following result:
   \begin{theorem}\label{Existence1}
   Let $y_0,v_0\in \R$ and $r>\max\{|y_0|,|v_0|\}$.
   Let $(\Phi_2)$, $(A_0)$ and $(A_l)$ hold.
   Then there exists $\delta=\del(r) > 0$,
   such that:
   \begin{itemize}
   \item[a)] $(CP)$ has a non-impact solution in $[0,\del]$.
   \item[b)] If $(A_l)'$ holds and $y_0=0$,
   then $(CP)$ has a classical solution in $[0,\del]$.
   \end{itemize}
  \end{theorem}

  \pr
   We consider operator $T$ defined by (\ref{opT}) and
   proceed in the same way as in Theorem \ref{Existence}.
   The proof of this theorem is very similar to the proof of Theorem
   \ref{Existence}, therefore we will just give
   the parts of the proof which are different.\\
  %\pagebreak

  \emph{(i) $T$ maps $C^i([0,a])$ into $C^i([0,a])$ for $i=1,2$.}\\

  Let $y\in C^1([0,a])$. The continuity of $Ty$
  follows as in Theorem
  \ref{Existence}. Let us  show the continuity
  of its first derivative $\frac{d}{dt}Ty$.

  An expression for operator $\frac{d}{dt}J$ (see  (\ref{prviizvodJ}))
  has an  additional summand, which does not vanish (we do not have
  $y(0)=0$),
  but continuity follows from the fact that
  $\phi_1(\ga)=0$ for $\ga \geq 1$.
  Therefore, %$\frac{d}{dt}T$:
  \begin{align}\label{Tizvod1}
  \frac{d}{dt}Ty(t) & =  \frac{d}{dt}(l*Jy)(t) + Gy(t)+ v_0,\\
   &= l*Jy'(t) + l*y(0)\int_0^d\frac{\phi_1(\ga)}{\Ga(2-\ga)}t^{1-\ga}d\ga +
                    \int_0^t f(t,y(t))+ v_0,\quad t\in [0,a],\nonumber
  %& = \frac{d}{dt}Ty(t) + y(0) \int_0^d\frac{\phi_1(\ga)}{\Ga(2-\ga)}l(t)*t^{1-\ga}d\ga.
  \end{align}
  and thus $T:  C^1([0,a])\to C^1([0,a])$.

  Let $y\in C^2_0([0,a])$ and assume $(A_l)'$ instead of $(A_l)$.
  Then, $Ty(0)=0$ because $y(0)=0$  and the first derivative of $Ty$ is given
  by (\ref{prviizvodT}).  The second derivative %$\frac{d^2}{dt^2}Ty$ is
  given by
  \begin{align}%\label{Tizvod2}
  \frac{d^2}{dt^2}Ty(t) & = l*Jy''(t) +
   l*y'(0)\int_c^d\frac{\phi_1(\ga)}{\Ga(2-\ga)}t^{1-\ga}d\ga
  + f(t,y(t)), \quad t\in [0,a]\nonumber
  % & = \frac{d^2}{dt^2}Ty(t)+ y'(0)\int_0^d\frac{\phi_1(\ga)}{\Ga(2-\ga)}l(t)*t^{1-\ga}d\ga
  \end{align}
  is continuous since $y''\in C([0,a])$, $J: C([0,a])\to C([0,a])$,
  $\phi_1 \equiv 0$ for $\ga \geq 1$ and $f$ is continuous.
  Therefore $ T:C^2_0([0,a])\to C^2_0([0,a])$.\\
%%%%%%%%%%%%%%%%%%%%%%%%%%%%%%%%%%%%%%%%%%%%%%%%%%

 \emph{(ii) {Additional regularity properties of $T$.}}\\

 Additional properties of $T$ which assure higher regularity of the
 solutions are to
  follow. We will show that $(A_l)$ implies the mapping properties
  \begin{equation}\label{absT}
  T: C([0,a])\to C^1([0,a])\quad\text{ and }\quad T: C^1([0,a])\to AC^2([0,a]).
  \end{equation}
  We will use the  following lemma.
  \begin{lemma}\label{absconofIal}
   a) Let $y\in C([0,a])$ and $\al \geq 1$. Then $I^{\al}y\in C^1([0,a])$.\\
   b) Let $y\in L^1_{\rm loc}([0,a])$ and $\al \geq 1$. Then $I^{\al}y\in AC([0,a])$.
  \end{lemma}
   \pr a)
   If $\al = 1$ then $\frac{d}{dt} I^{\al}y(t) = y(t)$ and therefore if $y\in C([0,a])$ then
   $I^{\al}y\in C^1([0,a])$, and if $y\in L^1_{\rm loc}([0,a])$
   then $I^{\al}y\in AC([0,a])$.
   For $\al>1$ using that $I^{\al}I^{\be}=I^{\al+\be}$ for $\al,\be > 0$
   and  $D^1I^1=I$
    we have that
    $\frac{d}{dt} I^{\al}y(t) = D^1 I^1 I^{\al-1}y(t)=I^{\al-1}y(t)$.
    If $y\in C([0,a])$ then $I^{\al-1}y\in C([0,a])$,
    and again $I^{\al}y\in C^1([0,a])$.
    If $y\in L^1_{\rm loc}([0,a])$ then $I^{\al-1}y \in
    L^1_{\rm loc}([0,a])$ (since $I^{\al-1}:L^1_{\rm loc}([0,a])\to
    L^1_{\rm loc}([0,a])$) and therefore $I^{\al}y\in AC([0,a])$.
   \ep

   Let $y\in C([0,a])$. Lemma \ref{absconofIal} gives $I^{\al}y\in C^1([0,a])$ if $\al > 1$,
  hence
  \begin{equation}\label{absJ}
   J:C([0,a])\to C^1([0,a]).
  \end{equation}
  Since $Jy(t) = \int_c^d\phi_1(\ga)D^{\ga-2}y\ d\ga, t\in [0,a] $ and $d<1,$ it follows that
   operator $D^{\ga-2}$ is in fact $I^{2-\ga}$ with $2-\ga > 1$; so  $D^{\ga-2}y\in C^1([0,a])$.
   Since $Gy\in C^1([0,a])$ we have that $Ty\in C^1([0,a])$.

  Let $y\in C^1([0,a])$.
  First note that if $\ga<1$ then $(t^{1-\ga})' = (1-\ga)t^{-\ga}\in
  L^1_{\text{loc}}(\R)$, hence
  \begin{equation}\label{Jdodatak}
  y(0)\int_c^d\frac{\phi_1(\ga)}{\Ga(2-\ga)}t^{1-\ga}d\ga\in
  AC([0,a]).
  \end{equation}
  Assumption $y'\in C([0,a])$, as well as  (\ref{absJ}) and (\ref{Jdodatak}) imply that
   $$\frac{d}{dt}Jy = Jy' + y(0)\int_c^d\frac{\phi_1(\ga)}{\Ga(2-\ga)}t^{1-\ga}d\ga
   \in AC([0,a]).$$
  Thus, $
  J: C^1([0,a])\to AC^2([0,a]).
  $
  Further, $Gy\in AC^2([0,a])$ because
  $\frac{d^2}{dt^2}Gy(t)=f(t,y(t))$, $t\in[0,a]$ and for $y$
  continuous  $t\mapsto f(t,y(t))$ is a locally integrable function.
  Therefore,  by (\ref{Tizvod1}) we have that
  %\begin{equation}\label{absT1}
  $T:  C^1([0,a])\to AC^2([0,a])$.\\
  %\end{equation}.

  %%%%%%%%%%%%%%%%%%%%%%%%%%%%%%%%%%%%%%%%%%%%

  \emph{(iii) $T$ is compact in $C^i([0,a])$, $i=1,2.$}\\

  Let $i=1.$ %The compactness of $T$ in $C^1([0,a])$,
 By Theorem \ref{Existence},
 $T$, given by (\ref{opT}) (with (\ref{prviizvodT})),
 maps a bounded sequence in $C^1([0,a])$ into
 a sequence with a convergent subsequence in $C^1([0,a])$.

 For the compactness of the operator $T$, given by (\ref{opT})
 (with (\ref{Tizvod1})),
 we need additionally to show  that if $\{y_n\}_{n\in\N}$
 is bounded in $C^1([0,a])$,
 then $$\{y_n(0)
 \int_c^d\frac{\phi_1(\ga)}{\Ga(2-\ga)}l(t)*t^{1-\ga}d\ga\}_{n\in\N}$$
 has  a convergent subsequence.
 This is true since any bounded sequence $\{y_n\}_{n\in\N}$ in $C^1([0,a])$ has a subsequence,
 again denoted
 by $\{y_n\}_{n\in\N}$, which is convergent in $C([0,a])$.
 Hence, $\{y_n(0)\}_{n\in\N}$ converges as well.

 Similarly, one can prove the compactness of $T$ in
 $C^2_0([0,a])$.\\

 \emph{(iv) Determination of $\del.$}\\

 Take  $B_r \subset C^1([0,\delta])$ and $B_r\subset C_0^2([0,\delta])$,
 respectively,
 and $\del$ will be defined later.
 Then, for $i\in \{0,1\}$,
 $$
 |l*y^{(i)}(0)\int_c^d\frac{\phi_1(\ga)}{\Ga(2-\ga)}t^{1-\ga}d\ga|
 \leq |y^{(i)}(0)|M'_t \leq  M'_t  \|y\|_{i+1},  \quad t\in
 [0,\del],
 $$
 where
 $$ M'_t = \int_0^t l(\tau)d\tau
 \int_c^d\frac{|\phi_1(\ga)|}{|\Ga(2-\ga)|}t^{1-\ga}d\ga,\quad t\in
 [0,\delta].$$
 Using this, (\ref{Tizvod1}) and estimate (\ref{ocenaT}),
 similarly as in the proof of Theorem \ref{Existence}, one can estimate
 $\frac{d^i}{dt^i}Ty(t)$, for $i=0,1$, if $(A_l)$ is assumed,
 and for $i=0,1,2$, if $(A_l)'$
 is assumed.
 Therefore, it is possible to find $\delta>0$ such that $y\in B_r\subset C^1([0,\del])$
 and  $y\in B_r\subset C^2_0([0,\del]) $, respectively,
 implies that
 $$
 |\frac{d^i}{dt^i}Ty(t)|\leq r,\quad i=0,1, \text{ respectively } i=0,1,2 .
 $$

%\pagebreak

\emph{ (v) Assertions a) and b).}\\

 a) The Schauder fixed point theorem implies the existence of
 a fixed point $y$ for $T$ in $C^1([0,\delta])$, if $(A_l)$ is
 assumed.  Let $z$ be given by (\ref{z1}). Then $(y,z)$ is mild
 solution to $(CP)$.
 By (\ref{absT}) we obtain that $y\in AC^2([0,\del])$. This, with $d<1$,  implies
 \begin{equation*}%\label{D2J}
   \int_c^d \phi_1(\ga)D^{\ga}y(t) d\ga\in AC([0,\del]).
  \end{equation*}
  Since $l$ is  locally integrable, % by (\ref{z1})
  we have that  $z\in AC([0,\del])$. Therefore, the mild solution $(y,z)$
  satisfies the first equation in $(CP)$ in $L^1_{\rm
  loc}([0,\del])$, so this is also a non-impact solution
  %Moreover,
  %if $\ga_0 < 1$ then the second equation in $(CP)$ also holds in $L^1_{\rm loc}([0,\del])$.

 b) If we assume $(A_l)'$, then $T$ has a fixed point in $C^2_0([0,\delta])$,
  $(CP)$ has a mild solution $(y,z)$ in $[0,\del]$, $z$ is given by
  (\ref{z1}),
  and $y\in C^2_0([0,\delta])$. So $(\widetilde{CP})_2$ holds.
  Further, $d<1$ implies that
  %\begin{equation}\label{D2J1}
  $\int_c^d\phi_1(\ga)D^{\ga}y(t) d\ga \in C^1([0,\del])$
  %\end{equation}
   and therefore, $z\in C^1([0,\del])$ and $(\widetilde{CP})_1$ is satisfied for
  all $t\in [0,\del]$.
   \ep

%%%%%%%%%%%%%%%%%%%%%%%%%%%%%%
%%%%%%%%%%%%%%%%%%%%%%%%%%%
 In this particular case, one can also impose
 a condition on $f$, stronger then $(A_l)$ and $(A_l)'$, and obtain
 classical solution to $(CP)$ in $[0,\del]$, which satisfies
 $(\widetilde{CP})_2$ for all $t\in [0,\del]$. The condition reads:
 \begin{align*}
 (A_l)''\quad & f(t,u),\,(t,u)\in [0,\infty )\times \R, \text{ is continuous }
  f(t,0)=0, t\in[0,\infty),\text{ and} \\
  & \text{ there exist } a>0,\, h\in C([0,a])\text{ with }h\geq 0
    \text{ and }h(0,0)=0,
  \text{ such that } \\
  & |f(t,u)-f(s,v)|\leq h(t,s)(|t-s| + |u-v|),\,u,v\in \R,t\in [0,a].
 \end{align*}
 The corresponding theorem is the following.
 \begin{theorem}\label{Existence11}
   Let $v_0\in \R$, $r> |v_0|$ and $y_0=0$.
   Let $(\Phi_2)$, $(A_0)$ and $(A_l)''$ hold.
   Then there exists $\delta=\del(r) > 0$,
   such that there exists the classical solution to $(CP)$
   in $[0,\del]$, which satisfies equation $(CP)_2$ point-wisely for $t\in [0,\del]$.
 \end{theorem}
 \pr
  First note that all assumptions of Theorem \ref{Existence1} are
  satisfied. Hence there is a classical solution to $(CP)$ in $[0,\del]$ obtained in
  five steps in previous proof.
  The assumption  $(A_l)''$ gives stronger results. Note that
  if $y\in AC([0,a])$ and $(A_l)''$ holds, then
  $t\mapsto f(t,y(t))$, $t\in [0,a]$ is absolutely continuous.
  Indeed,
  \begin{align*}
   |f(t,y(t))-f(s,y(s))|&\leq h(t,s)(|t-s| + |y(t)-y(s)|)
  \end{align*}
  implies absolutely continuity since  $y$ belongs to $AC([0,a])$.

  Further note that in step \emph{(ii)} of the proof of the Theorem
  \ref{Existence1}, in addition to (\ref{absT})
  we have
  \begin{equation}\label{absT2}
   T: C^2([0,a])\to AC^3([0,a]).
  \end{equation}
  To see this, let $y\in C^2_0([0,a])$. Then (\ref{absJ}) and (\ref{Jdodatak}) imply that
  $$\frac{d^2}{dt^2}Jy = Jy'' +  y'(0)\int_c^d\frac{\phi_1(\ga)}{\Ga(2-\ga)}t^{1-\ga}
  \in AC^3([0,a]).$$
  Since $t \mapsto f(t,y(t))$, $t\in[0,a]$ is absolutely
  continuous, we have that  $Gy\in AC^3([0,a])$ and
  therefore (\ref{absT2}). It follows that the fixed point $y \in C_0^2([0,a])$
 obtained in part b) of the fifth step of the proof of the Theorem \ref{Existence1}
 is then an element of $AC^3([0,\del])$.
 Then $z$ given by (\ref{z1}) is an element in $AC^2([0,\del])$,
 which yields that both equations in $(CP)$ are satisfied for all
 $t\in [0,\del]$.
 \ep

 \begin{remark}\rm
 a) Note that in Theorem \ref{Existence} and
  Theorem \ref{Existence1} we could have assumed
  that $[c,d]$ was a subset of the interval $(-\infty,2)$ and
  $(-\infty,1)$, respectively,
  in which case we would have the same results.
  Indeed, with respective assumptions, we change the bounds
  of the integral in (\ref{J}), the  definition  of operator $J$,
  but not its properties: $\frac{d^i}{dt^i}Jy$, $i=0,1,2$
  remain continuous, $J$ remains compact in respective spaces, and
  estimates in (\ref{ocenaT}) hold.

 b) With assumptions $d<2$ and $d<1$, respectively we
   could also have  considered a constitutive equation of the form
  $$
  \sum_{i=0}^k a_iD^{\ga_i}z(t)=\int_{-\infty}^d \phi_1(\ga)D^{\ga}y(t)d\ga
  $$
  and we would obtained the same results as in Theorem \ref{Existence}
  and in Theorem  \ref{Existence1},
  respectively.

  c) It also make a sense to consider constitutive equation for
  $d<0$. Then  on its right hand side only
  fractional integrals  of $y$ appears. Such case is covered with
   condition that follows.
 \end{remark}

  \begin{align*}
  (\Phi_3) \quad
  &\phi_1 \text{ is  continuous function in $[c,d]\subset (-\infty,0)$,
  $\phi_1(\ga)=0$, $\ga\notin [c,d]$};\\
  &\phi _2 \text{ is of the form (\ref{fidva})}.
  \end{align*}
  Also,  the theorem similar to Theorems \ref{Existence} and
  \ref{Existence1} holds.
  \begin{theorem}\label{Existence2}
   Let $y_0,v_0\in\R$ and $r>\max\{|y_0|,|v_0|\}$.
   Let $(\Phi_3)$, $(A_0)$  and $(A_l)'$ hold.
   Then there exists $\delta=\del(r) > 0$ such that
   $(CP)$ has a classical solution in $[0,\del]$.
   Moreover, the classical solution satisfies
   $(CP)_2$ for all $t\in [0,\del]$.
  \end{theorem}
 %\subsection{$(C_{\phi_1,\phi_2})^4 $ - case}

 Finally, we are interested in cases when  both
 $\phi_1$ and $\phi_2$ are linear combinations of translations of delta
 distributions:
 \begin{align*}
 (\Phi_4) \quad
 &\phi_1(\cdot) = \sum_{j=0}^m b_j\del(\cdot-\be_j),\be_i\in [0,2),
  \be_0\geq \be_j \geq \be_m,  j\in\{0,...,m\}\\
 &\phi _2 \text{ is of the form (\ref{fidva})}.
  \end{align*}
 The constitutive equation $(CP)_2$ becomes
 $$
  \sum_{i=0}^k a_iD^{\ga_i}z(t)=\sum_{j=0}^m b_jD^{\be_j}y(t),\quad
  t>0,
 $$
 and the theorem similar to previous holds.
   \begin{theorem}\label{Existence3}
  Let $y_0,v_0\in\R$ and $r>\max\{|y_0|,|v_0|\}$.
  Let $(\Phi_4)$, $(A_0)$ and $(A_l)$ hold.
  Then there exists $\del=\del(r) > 0$ such that:
 \begin{itemize}
   \item[a)] $(CP)$ has a mild solution in $[0,\del]$.
   \item[b)] $(CP)$ has a non-impact solution in $[0,\del]$ if one of the following
    conditions hold
   $$(i)\, y_0=0 \text{ and } \ga_0-\ga_k>1 \quad \text{ or }\quad (ii)\, \be_0 < 1.$$
   \item[c)] If $(A_l)'$ holds, then $(CP)$ has
    a classical solution in $[0,\del]$ if one of the following
    conditions hold
    $$ (i)\, y_0 = v_0 = 0\quad \text{ or }\quad (ii)\, \be_0 < 1 \text{ and } y_0=0
    \quad  \text{ or }\quad (iii)\, \be_0 < 0.$$
 \end{itemize}
 \end{theorem}
 The proof similar to previous ones is skipped.

 \section{Continuous $\phi_2$}\label{SC2}
 % \todo{ili naslov Case $\phi_2\in C^3([0,1])$}

 In this section we impose different condition to $\phi_2$.

% To show the existence of the solutions in this case we consider a more general one:

  \begin{align*}
 (\Phi_5): \quad
 \bullet\ &\phi_1 \text{ is continuous function in } [c,d]\subset [0,2),\,
 \phi_1(\ga)=0,\, \ga\notin [c,d]\text{ or }\\
 &\phi_1(\cdot) = \sum_{j=0}^m b_j\del(\cdot-\be_j),\be_i\in [0,2),
  \be_0\geq \be_j \geq \be_m,  j\in\{0,...,m\}\\
 \bullet\ &\phi_2\in C^3([0,1]),\, \phi_2\equiv 0 \text{ out of }
 [0,1],\, \phi_2(1)\neq 0 \text{ and}\\
 & \text{ either }\phi(0)\neq 0 \text{ or } \phi_2(\ga)\sim p\ga^q,\, p>0,\, q>0.
  \end{align*}

 \noindent Assumptions on $\phi_2$ given in  $(\Phi_5)$
 will imply (as it is shown in \cite{Kochubei:07}) the existence of a solution to
 distributed order differential equation $\int_0^1 \phi_2(\ga)D^{\ga} y = g$, i.e.
 the existence of a locally integrable function
 $\ILT(\frac{1}{\int_0^1\phi_2(\ga)s^{\ga}d\ga})$.

 \begin{theorem}\label{Existence4}
  Let $y_0, v_0\in\R$, $r>\max\{|y_0|,|v_0|\}$ and $f$ satisfy $(A_l)$.
  Let $\phi_1$ and $\phi_2$ satisfy $(\Phi_5)$.
  Then there exists $\del=\del(r)>0$  such that:\\
  a) There exists a mild solution of $(CP)$ in $[0,\del]$.\\
  b) There exists a non-impact solution of $(CP)$ in  $[0,\del]$
   if $\supp \phi_1\subset [c, 1)$.\\
  c) There exists a classical solution of $(CP)$ in
   $[0,\del]$ if $(A_l)'$  and one of the following conditions
   hold:\\
  $ (i)\  y_0=v_0=0\,\text{  or   }\, (ii)\ \supp \phi_1\subset [c, 1)
    \text{   and  } y_0=0\,  \text{   or   }\,  (iii)\ \supp \phi_1\subset [c,
    0)$.
 \end{theorem}
 \pr
 The solution to equation
 \begin{equation*}
  \int_0^1 \phi_2(\ga)D^{\ga}zd\ga = g,\quad \text{ in }\spp
 \end{equation*}
 is given by
 $$
 z = \chi*g, \text{ where }
 \chi=\ILT(\frac{1}{\int_0^1\phi_2(\ga)s^{\ga}d\ga}),
 $$
 provided that the inversion exist.
 It is proved in \cite[Proposition 3.1]{Kochubei:07},
 that conditions on $\phi_2$ given in $(\Phi_5)$,
  imply that
 $\chi$ exists and,
 moreover, that $\chi\in C^{\infty}((0,\infty))\cap L^1_{\text{loc}}([0,\infty])$.
 Therefore, $(CP)_2$
 $$
 \int_0^1 \phi_2(\ga)D^{\ga}zd\ga = \int_0^2
 \phi_1(\ga)D^{\ga}yd\ga
 $$
 has a solution in $\spp$
 $$
 z = \chi*\int_0^2 \phi_1(\ga)D^{\ga}yd\ga.
 $$
 The substitution of $z$ in $(CP)_1$ and the integrations give
 \begin{equation*}
 y = -\chi*\int_0^2\phi(\ga)D^{\ga-2}yd\ga + D^{-2}f(\cdot,y(\cdot))
 + v_0 t + y_0,\quad \text{ in }\spp.
 \end{equation*}
 Again, we consider operator $T$ which is given by (\ref{opT}),
 where $R$ is now given by
 $
 Ry := \chi*Jy
 $
 and $J$ by (\ref{J}). Since $\chi$ is locally
 integrable, the same properties hold for $J$, $R$ and $T$, as in Theorems
 \ref{Existence} and others.
 Thus, it follows that there exist $\del=\del(r)$
 and a mild solution $(y,z)$ to $(CP)$ in $[0,\del]$.
 Also, one can prove the properties of  $y$ and $z$
 quoted  in a), b), and c)  in the same way as in Theorems
 \ref{Existence} and others.\\
 \ep

% \noindent \textbf{Example: $\phi_1 = a^{\ga}$, $\phi_2= b^{\ga}$}\\

 To close the section we give the example that follows arise from application. Consider the
 system $(CP)$ with
 $\phi_1(\ga) = b^{\gamma}$ and $\phi_2(\ga) = a^{\gamma}$, $\ga\in
 [0,2)$, where $a$ and $b$  are positive constants with $b > a$. The
 latter condition is consequence of the Second  Law of
 thermodynamics.
 \begin{equation*}
  (CP)_{ex}\quad\quad
  \left\{ \begin{array}{c}
   D^{2}y(t)+z(t)=f(t,y(t)),\quad t > 0, \\
   \\
   \int_0^1 b^{\ga}D^{\ga}y(t)d\ga = \int_0^1
   a^{\ga}D^{\ga}z(t)d\ga, \quad t>0\\
   \\
   y(0)=y_{0};\qquad y^{(1)}(0)=v_{0},
  \end{array}
  \right.
 \end{equation*}
 One can apply Theorem \ref{Existence4} to prove the existence of
 a mild and a classical solution for above system.
 \begin{theorem}\label{existenceEX}
  Let $y_0, v_0\in\R$ and $r>\max\{|y_0|,|v_0|\}$. Let $(A_l)$ hold.
  Then problem  $(CP)_{ex}$\\
     (i) has a non-impact solution.\\
   (ii) has a classical solution $(y,z)$ if  $(A_l)'$ holds and
     if  $y_0=0$.
   \end{theorem}
 \pr
  Since $\phi_1,\phi_2\in C^3([0,1])$, $\phi_2(1)=b\neq 0$
  and  $\phi_2(0)=1\neq 0$
  conditions of Theorem  \ref{Existence4} (part b) case (ii)  and part c) case
  (ii))
  are satisfied. Theorem \ref{Existence4} implies the assertions.
 \ep
% \bibliography{ljubica}
 %\bibliographystyle{abbrv}
 % \newcommand{\SortNoop}[1]{}
\newcommand{\SortNoop}[1]{}

 \end{document}